\documentclass[12pt]{amsart}
\usepackage{amssymb}

\expandafter\ifx\csname delta.sty\endcsname\relax \else\endinput\fi
\expandafter\edef\csname delta.sty\endcsname{%
 \catcode`\noexpand\@=\the\catcode`\@\space}

\let\atbefore @

\catcode`\@=11

\newif\ifMag
\let\@ft@\expandafter
\numberwithin{equation}{section}

\newif\ifComments

\def\tod@y{\ifcase\month\or
 January\or February\or March\or April\or May\or June\or July\or
 August\or September\or October\or November\or December\fi\space\,
\number\day,\space\,\number\year}

\def\h@@r{hh}\def\m@n@te{mm}
\def\wh@tt@me{\count@\time\divide\count@ 60\edef\h@@r{\number\count@}%
 \multiply\count@ -60\advance\count@\time\edef
 \m@n@te{\ifnum\count@<10 0\fi\number\count@}}
\def\t@me{\h@@r\/{\rm:}\m@n@te} \let\whattime\wh@tt@me
\let\Today\tod@y \let\nowtime\t@me

\def\ftext#1{{\let\thefootnote\relax\footnotetext{\vsk-.8>\nt #1}}}
\def\em#1{{\itshape #1\/}}

 {\end{trivlist}\egroup\global\@ignoretrue\ifComments\unvbox\commentbox\fi
 \noindent}

\def\gadv{\global\adv} \def\gad#1{\gadv#1\@ne} \def\gadneg#1{\gadv#1-\@ne}
\def\textindent#1{\indent\llap{#1\enspace}\ignorespaces}

\newcount\itemlet
\def\newbi{\itemlet 96} \newbi
\def\bitem{\gad\itemlet\endgraf\hangindent1.5\parindent
 \hglue-.5\parindent\textindent{\upshape\rlap{\char\the\itemlet}\hp{b})}}

\newcount\itemrm

\def\iitem{\gad\itemrm\endgraf\hangindent1.5\parindent\hglue-.5\parindent
 \textindent{\upshape\hp{v}\llap{\romannumeral\the\itemrm})}}

\let\Disp\[ \let\endD\] \let\{\protect

\def\Tag#1{\label{e:#1}\let\notag\relax} 
\def\sh@nd#1#2{\begin{#1*}#2\end{#1*}}
\def\n@t@gs#1#2#3{\let\n@@@l\\ \begin{#1}#2\global\let\d@bl\\
 \gdef\\{\notag\d@bl}#3\notag\global\let\\\d@bl\end{#1}\let\\\n@@@l}

\def\Gather#1\endG{\n@t@gs{gather}{}{#1}}
\def\gAther#1\endG{\sh@nd{gather}{#1}}
\def\Align#1\endA{\n@t@gs{align}{}{#1}}
\def\aLign#1\endA{\sh@nd{align}{#1}}
\def\Alignat#1#2\endAt{\n@t@gs{alignat}{#1}{#2}}
\def\aLignat#1\endAt{\sh@nd{alignat}{#1}}

\def\(#1){\textup{(\ref{e:#1})}}
\def\[{\@ifnextchar:\c@t@sect\c@t@}
\def\c@t@sect:#1]{\ref{s:#1}} \def\c@t@#1]{\ref{t:#1}}

\def\qed{\hbox{}\nobreak\hfill\nobreak{\m@th$\,\square$}}
\def\sk@@p#1{\par\skip@#1\relax\ifdim\lastskip<\skip@\relax\removelastskip
 \vskip\skip@\fi}
\def\demo#1{\sk@@p\medskipamount\nt{\ignore\it #1\unskip.}\enspace
 \ignore}
\def\enddemo{\sk@@p\medskipamount}
\def\Rem{\demo{\sl Remark}} 
\def\Pf#1.{\demo{Proof #1}}

 \let\nl\newline
\let\bls\baselineskip \let\ignore\ignorespaces \let\adv\advance
\def\vsk#1>{\vskip#1\bls}
\def\vv#1>{\vadjust{\vsk#1>}\ignore}
\def\vvn#1>{\vadjust{\nobreak\vsk#1>\nobreak}\ignore}
\def\vvv#1>{\vskip\z@\vsk#1>\nt\ignore}
\def\vvgood{\vadjust{\penalty-500}} 
\def\gooddb{\noalign{\penalty-500}}
\def\mathbox#1{\hbox{\m@th$#1$}}

\let\dsize\displaystyle \let\tsize\textstyle
\let\ssize\scriptstyle \let\sss\scriptscriptstyle
 
\let\vp\vphantom \let\hp\hphantom \let\nt\noindent
\def\hline{\hbox to\hsize}
\let\cline\centerline \let\lline\leftline \let\rline\rightline
\def\nn#1>{\noalign{\vskip#1\p@}} \def\NN#1>{\openup#1\p@}
 
\let\Lim\lim \def\lim{\Lim\limits} \let\Sum\sum \def\sum{\Sum\limits}
 
\let\Prod\prod \def\prod{\Prod\limits} \let\Int\int \def\int{\Int\limits}

\def\tsum{\mathop{\tsize\Sum}\limits} 

\def\~{\leavevmode\@ifnextchar~\m@n@s\@md@sh}
\def\m@n@s~{\raise.15ex\mathbox{-}} \def\@md@sh{\raise.13ex\hbox{--}}
\let\procent\% \def\%#1{\ifmmode\mathop{#1}\limits\else\procent#1\fi}
\let\@ml@t\" \def\"#1{\ifmmode ^{(#1)}\else\@ml@t#1\fi}
\let\@c@t@\' \def\'#1{\ifmmode _{(#1)}\else\@c@t@#1\fi}
\let\colon\: \def\:{^{\vp|}} \def\&{.\kern.1em} \def\^#1{\text{\m@th#1}}

\newif\ifNewskips

\def\Newskips{\global\Newskipstrue
 \gdef\>{\relax\ifmmode\mskip.666667\thinmuskip\relax\else\kern.111111em\fi}
 \gdef\}{\relax\ifmmode\mskip-.666667\thinmuskip\relax\else\kern-.111111em\fi}
 \gdef\){\relax\ifmmode\mskip.333333\thinmuskip\relax\else\kern.0555556em\fi}
 \gdef\]{\relax\ifmmode\mskip-.333333\thinmuskip\relax\else\kern-.0555556em\fi}}
\Newskips

\def\Re{\mathop{\mathrm{Re}\>}} \def\Im{\mathop{\mathrm{Im}\>}}
\def\End{\mathop{\mathrm{End}\>}} 
 \def\Sym{\mathop{\mathrm{Sym}\)}}

\def\1{^{-1}} \def\_#1{_{\Rlap{#1}}}
\def\vst#1{{\lower1.9\p@
 \hbox{\m@th$\bigr|_{\raise.5\p@\hbox{\m@th$\ssize#1$}}$}}}
\def\vrp#1:#2>{{\vrule height#1 depth#2 width\z@}}
\def\vru#1>{\vrp#1:\z@>} \def\vrd#1>{\vrp\z@:#1>}
 
\def\sscr#1{\raise.3ex\hbox{\m@th$\sss#1$}} \def\@@PS{\mathbf{OOPS!!!}}

\def\lsym#1{#1\alb\ldots\relax#1\alb}
\def\lc{\lsym,}   \def\lox{\lsym\ox}

\let\texspace\ \def\ {\ifmmode\alb\fi\texspace}

\def\Line#1{\kern-.5\hsize\hline{\m@th$\dsize#1$}\kern-.5\hsize}
\def\Lline#1{\kern-.5\hsize\lline{\m@th$\dsize#1$}\kern-.5\hsize}
\def\Cline#1{\kern-.5\hsize\cline{\m@th$\dsize#1$}\kern-.5\hsize}
\def\Rline#1{\kern-.5\hsize\rline{\m@th$\dsize#1$}\kern-.5\hsize}

\def\Ll@p#1{\llap{\m@th$#1$}} \def\Rl@p#1{\rlap{\m@th$#1$}}
 \def\Cl@p#1{\llap{\m@th$#1$\hss}}
\def\Llap#1{\mathchoice{\Ll@p{\dsize#1}}{\Ll@p{\tsize#1}}{\Ll@p{\ssize#1}}%
 {\Ll@p{\sss#1}}}
\def\Clap#1{\mathchoice{\Cl@p{\dsize#1}}{\Cl@p{\tsize#1}}{\Cl@p{\ssize#1}}%
 {\Cl@p{\sss#1}}}
\def\Rlap#1{\mathchoice{\Rl@p{\dsize#1}}{\Rl@p{\tsize#1}}{\Rl@p{\ssize#1}}%
 {\Rl@p{\sss#1}}}
 
\def\LRtph#1#2{\setbox\z@\hbox{#1}\dimen\z@\wd\z@\hbox{\hbox to\dimen\z@{#2}}}
\def\LRph#1#2{\LRtph{\m@th$#1$}{\m@th$#2$}}

\def\Lto#1{\setbox\z@\hbox{\m@th$\tsize{#1}$}%
 \mathrel{\mathop{\hbox to\wd\z@{\rightarrowfill}}\limits#1}}
\def\Lgets#1{\setbox\z@\hbox{\m@th$\tsize{#1}$}%
 \mathrel{\mathop{\hbox to\wd\z@{\leftarrowfill}}\limits#1}}

\def\vpb#1{{\vp{\big(}}^{\]#1}}

\let\alb\allowbreak

 \let\x\times \let\ox\otimes 
\let\sub\subset 
\let\le\leqslant \let\ge\geqslant
 \let\8\infty \let\*\star

 \def\vert{\ |\ }

\let\lb\lbrace \let\rb\rbrace

\let\Bbb\mathbb

\let\Cal\mathcal
\let\frak\mathfrak

\def\pms{\raise.25ex\mathbox{\ssize\pm}\>}
\def\mps{\raise.25ex\mathbox{\ssize\mp}\>}

\let\al\alpha
\let\bt\beta
\let\gm\gamma \let\Gm\Gamma 
\let\dl\delta  
 \let\eps\varepsilon \let\epsilon\eps

\let\ka\kappa
\let\la\lambda

 \let\phi\varphi

\def\C{\Bbb C}
\def\R{\Bbb R}
\def\Z{\Bbb Z}

\def\Zp{\Z_{\ge 0}}

\def\h@ph{\discretionary{}{}{-}} \def\$#1$-{\,\text{\m@th$#1$}\h@ph}

\def\difl/{differential} \def\dif/{difference}
\def\cf.{cf.\ \ignore} \def\Cf.{Cf.\ \ignore}
\def\egv/{eigenvector} \def\eva/{eigenvalue} \def\eq/{equation}
\def\lhs/{the left hand side} \def\rhs/{the right hand side}
\def\Lhs/{The left hand side} \def\Rhs/{The right hand side}
\def\gby/{generated by} \def\wrt/{with respect to} \def\st/{such that}
\def\resp/{respectively} \def\off/{offdiagonal} \def\wt/{weight}
\def\pol/{polynomial} \def\rat/{rational} \def\tri/{trigonometric}
\def\fn/{function} \def\var/{variable} \def\raf/{\rat/ \fn/}
\def\inv/{invariant} \def\hol/{holomorphic} \def\hof/{\hol/ \fn/}
\def\mer/{meromorphic} \def\mef/{\mer/ \fn/} \def\mult/{multiplicity}
\def\sym/{symmetric} \def\perm/{permutation}
\def\rep/{representation} \def\irr/{irreducible} \def\irrep/{\irr/ \rep/}
\def\hom/{homomorphism} \def\aut/{automorphism} \def\iso/{isomorphism}
\def\lex/{lexicographical} \def\as/{asymptotic} \def\asex/{\as/ expansion}
\def\ndeg/{nondegenerate} \def\neib/{neighbourhood} \def\deq/{\dif/ \eq/}
\def\hw/{highest \wt/} \def\gv/{generating vector} \def\eqv/{equivalent}
\def\msd/{method of steepest descend} \def\pd/{pairwise distinct}
\def\wlg/{without loss of generality} \def\Wlg/{Without loss of generality}
\def\onedim/{one-dim\-en\-sion\-al} \def\fd/{fi\-ni\-te-dim\-en\-sion\-al}
\def\qcl/{quasiclassical} \def\hwv/{\hw/ vector}
\def\hgeom/{hyper\-geo\-met\-ric} \def\hint/{\hgeom/ integral}
\def\hwm/{\hw/ module} \def\emod/{evaluation module} \def\Vmod/{Verma module}
\def\symg/{\sym/ group} \def\sol/{solution} \def\eval/{evaluation}
\def\anf/{analytic \fn/} \def\anco/{analytic continuation}
\def\qg/{quantum group} \def\qaff/{quantum affine algebra}

\def\Rm/{\^{$R$-}matrix} \def\Rms/{\^{$R$-}matrices}
\def\YB/{Yang-Baxter \eq/}
\def\Ba/{Bethe ansatz} \def\Bv/{Bethe vector} \def\Bae/{\Ba/ \eq/}
\def\KZv/{Knizh\-nik-Zamo\-lod\-chi\-kov} \def\KZvB/{\KZv/-Bernard}
\def\KZ/{{\sl KZ\/}} \def\qKZ/{{\sl qKZ\/}}
\def\KZB/{{\sl KZB\/}} \def\qKZB/{{\sl qKZB\/}}
\def\qKZo/{\qKZ/ operator} \def\qKZc/{\qKZ/ connection}
\def\KZe/{\KZ/ \eq/} \def\qKZe/{\qKZ/ \eq/} \def\qKZBe/{\qKZB/ \eq/}

\def\LPT/{Laboratoire de Physique Th\'eorique ENSLAPP}
\def\ENSLyon/{\'Ecole Normale Sup\'erieure de Lyon}
\def\DMS/{Department of Mathematics, Faculty of Science}
\def\DMO/{\DMS/, Osaka University}
\def\DMOaddr/{Toyonaka, Osaka 560, Japan}
\def\dmoemail/{vt@math.sci.osaka-u.ac.jp}
\def\SPb/{St\&Petersburg}
\def\home/{\SPb/ Branch of Steklov Mathematical Institute}
\def\homeaddr/{Fontanka 27, \SPb/ \,191011, Russia}
\def\homemail/{vt@pdmi.ras.ru}
\def\absence/{On leave of absence from \home/}
\def\UNC/{Department of Mathematics, University of North Carolina}
\def\ChH/{Chapel Hill}
\def\UNCaddr/{\ChH/, NC 27599, USA} \def\avemail/{av@math.unc.edu}
\def\grant/{NSF grant DMS\~9501290}	
\def\Grant/{Supported in part by \grant/}

\def\Aomoto/{K\&Aomoto}
\def\Dri/{V\]\&G\&Drin\-feld}
\def\Fadd/{L\&D\&Fad\-deev}
\def\Feld/{G\&Felder}
\def\Fre/{I\&B\&Fren\-kel}
\def\Gustaf/{R\&A\&Gustafson}
\def\Kazh/{D\&Kazhdan} \def\Kir/{A\&N\&Kiril\-lov}
\def\Kor/{V\]\&E\&Kore\-pin}
\def\Lusz/{G\&Lusztig}
\def\MN/{M\&Naza\-rov}
\def\Resh/{N\&Reshe\-ti\-khin} \def\Reshy/{N\&\]Yu\&Reshe\-ti\-khin}
\def\Skl/{E\&K\&Sklya\-nin}
\def\SchV/{V\]\&\]V\]\&Schecht\-man} \def\Sch/{V\]\&Schecht\-man}
\def\Takh/{L\&A\&Takh\-tajan}
\def\VT/{V\]\&Ta\-ra\-sov} \def\VoT/{V\]\&O\&Ta\-ra\-sov}
\def\Varch/{A\&\]Var\-chenko} \def\Varn/{A\&N\&\]Var\-chenko}

\def\AMS/{Amer.\ Math.\ Society}
\def\CMP/{Comm.\ Math.\ Phys.{}}
\def\DMJ/{Duke.\ Math.\ J.{}}
\def\Inv/{Invent.\ Math.{}} 
\def\IMRN/{Int.\ Math.\ Res.\ Notices}
\def\JPA/{J.\ Phys.\ A{}}
\def\JSM/{J.\ Soviet\ Math.{}}
\def\LMP/{Lett.\ Math.\ Phys.{}}
\def\LMJ/{Leningrad Math.\ J.{}}
\def\LpMJ/{\SPb/ Math.\ J.{}}
\def\SIAM/{SIAM J.\ Math.\ Anal.{}}
\def\SMNS/{Selecta Math., New Series}
\def\TMP/{Theor.\ Math.\ Phys.{}}
\def\ZNS/{Zap.\ nauch.\ semin. LOMI}

\def\ASMP/{Advanced Series in Math.\ Phys.{}}

\def\AMSa/{AMS \publaddr Providence}
\def\Birk/{Birkh\"auser}
\def\CUP/{Cambridge University Press} \def\CUPa/{\CUP/ \publaddr Cambridge}
\def\Spri/{Springer-Verlag} \def\Spria/{\Spri/ \publaddr Berlin}
\def\WS/{World Scientific} \def\WSa/{\WS/ \publaddr Singapore}

\csname delta.sty\endcsname

\newif\ifUS
\ifx\USversion\relax\UStrue\fi

\ifUS
\fi

\evensidemargin\oddsidemargin
\edef\restoreatcode{\catcode`\noexpand\@=\the\catcode`\@}
\makeatletter

\newtheorem{theorem}{Theorem}
\newtheorem{proposition}[theorem]{Proposition}
\newtheorem{lemma}[theorem]{Lemma}
\newtheorem{corollary}[theorem]{Corollary}

\newenvironment{abstr}{\begingroup\narrower\small
\nt{\sc Abstract.}\enspace\ignorespaces}{\endgraf\endgroup}

\newenvironment{remark}{\demo{\sl Remark}}{\enddemo}
\newenvironment{example}{\demo{\bf Example}}{\enddemo}

\def\fratop{\genfrac{}{}{\z@}1}

\def\Ref#1{{\rm(\ref{#1})}}

\def\Ib{\bar I}
\def\Jb{\bar J}

\def\Cc{\Cal C}

\def\hg{\frak h}

\def\gl{\frak{gl}}
\def\gsl{\frak{sl}}
\def\glt{\gl_{\)2}}

\def\Et{\widetilde E}
\def\Ut{\widetilde U}

\def\emu{e^{\>\mu}}

\def\Gmf#1{\Gm\bigl(#1\bigr)}
\def\Gmk#1{\Gmf{(#1)/\ka}}
\def\Sin#1{\sin\)\bigl(#1\bigr)}
\def\sink#1{\Sin{\pi{(#1)}/\ka}}
\def\Exp#1{\exp\)\bigl(#1\bigr)}
\def\expk#1{\Exp{\pi i\)(#1)/\ka}}

\def\r-{\>\hbox{\m@th$r$-}}
\def\R-{\>\hbox{\m@th$R\)$-}}
\def\V-{\>\hbox{\m@th$V\}$-}}
\def\q-{\>\hbox{\m@th$q\)$-}}
\def\qtit-{\>\hbox{\protect\large\m@th$q$\)-}}
\def\qhead-{\>\hbox{\small\m@th$q$\)-}}

\restoreatcode

\title[\smash{Identities between \qhead-Hypergeometric and Hypergeometric
Integrals}] {Identities between \qtit-Hypergeometric and\\[3pt]
Hypergeometric Integrals of Different Dimensions}

\author[\smash{V\]\&Tarasov and A\&\]Varchenko}]
{\vbox{}V\]\&Tarasov$^\star$ \>and \;A\&\]Varchenko$^\diamond$}

\begin{document}

\hrule width0pt
\ifUS\else\vsk-.5>\fi

\maketitle

\ftext{\mathsurround 0pt
$\]^\star\)$Supported in part by RFFI grant 02\)\~\)01\~\)00085a
\>and \,CRDF grant RM1\~\)2334\)\~MO\)\~\)02\vv.2>\\
\hp{$^*$}{\normalsize\sl E-mail\/{\rm:} \homemail/}\vv.1>\\
${\]^\diamond\)}$Supported in part by NSF grant DMS\)\~\)0244579\vv.2>\\
\vv-1.2>
\hp{$^*$}{\normalsize\sl E-mail\/{\rm:} anv@email.unc.edu}}

\begin{center}
{\it\=$^\star$\home/\\[2pt]\homeaddr/\\[6pt]
$^\diamond$\UNC/ at \ChH/\\[2pt]\ChH/, NC 27599\~3250, USA}

\ifUS\vsk1.5>\else\vsk1.3>\fi
{\sl September 2003}
\end{center}

\ifUS\vsk1.2>\else\vsk>\fi

\begin{abstr}
Given complex numbers $m_1\),l_1$ and nonnegative integers $m_2\),l_2$,
such that $m_1+m_2=l_1+l_2$, for any $a\),b=0\lc\min\)(m_2\),l_2)$ we define
an \$l_2\)$-dimensional Barnes type \q-hypergeometric integral
$I_{a,b}(z,\mu\);m_1,m_2,l_1,l_2)$ and an \$l_2\)$-dimensional hypergeometric
integral $J_{a,b}(z,\mu\);m_1,m_2,l_1,l_2)$. The integrals depend on complex
parameters $z$ and $\mu$. We show that $I_{a,b}(z,\mu\);m_1,m_2,l_1,l_2)$
equals $J_{a,b}(\emu,z\);l_1,l_2,m_1,m_2)$ up to an explicit factor, thus
establishing an equality of \$l_2$-dimensional \q-hypergeometric and
\$m_2\)$-dimensional hypergeometric integrals. The identity is based on the
$(\gl_k\>,\gl_n)$ duality for the \qKZ/ and dynamical difference equations.
\end{abstr}

\thispagestyle{empty}

\ifUS\vsk>\else\vsk.5>\fi
\vsk0>
\section{Introduction}

\subsection{\qtit-Hypergeometric integrals}
\label{q-hyp}

Let $\ka$ be a positive number.
Let $m_1, l_1$ be complex numbers and $m_2, l_2$ nonnegative integers
such that
\vvn.2>
$$
m_1+\>m_2\,=\,l_1+\>l_2\,.
\vv.2>
$$
We say that an integer $a$ is \em{admissible} with respect to $m_2, l_2$ if
\vvn.2>
$$
0\le a\le\min\)(m_2, l_2)\,.
$$

\vsk.3>
For a pair of admissible numbers $a, b$ we define a function $I_{a,b}(z,\mu\);
m_1, m_2, l_1, l_2)$ of complex variables $z,\mu$. The function is defined as
an \$l_2\)$-dimensional Barnes type \q-hypergeometric integral:
\ifUS\vvn.1>\else\vvn.4>\fi
\begin{align}
\label{qmain}
& I_{a,b}(z,\mu\); m_1, m_2, l_1, l_2)\,={}
\\[6pt]
& {}\!=\!\int_{\dl_{l_2}(z\);m_1,m_2)}\!\!\!\!{}
\Phi_{l_2}(t,z,\mu\); m_1, m_2)\,
w_{l_2-a,a}(t,z\);m_1,m_2)\,W_{l_2-b,b}(t,z\);m_1,m_2)\,dt^{\>l_2}\,.
\notag
\end{align}
Here ${t=(t_1\lc t_{l_2})}$ and ${dt^{\>l_2}=dt_1\ldots\)dt_{l_2}}$.
\ifUS\vvgood\fi
The functions $\Phi_{l_2}(t,z,\mu\);m_1, m_2)$,
$w_{l_2-a,a}(t,\alb z\);\alb m_1, m_2)$ and $W_{l_2-b,b}(t, z\);m_1, m_2)$
are defined below. The \$l_2$-dimen\-sional integration contour
$\dl_{l_2}(z\);m_1,m_2)$ lies in $\C^{\>l_2}$ and is also defined below.
\vsk.2>
The \em{\q-master} function $\Phi_l$ is defined by the formula
\vvn.2>
\begin{align*}
& \Phi_l(t_1\lc t_l, z,\mu\); m_1, m_2)\,={}
\\[6pt]
& \!{}=\,\Exp{\)\mu\tsum_{u=1}^l\)t_u/\ka}\,
\prod_{u=1}^l\,\frac{\Gm(t_u/\ka)\,\Gmk{t_u\]-z}}
{\Gmk{t_u\]+m_1}\,\Gmk{t_u\]-z+m_2}}\,
\prod_{1\le u<v\le l}\,\frac{\Gmk{t_u\]-t_v\]+1}}{\Gmk{t_u\]-t_v\]-1}}\;.
\\[\ifUS-8pt\else-10pt\fi]
\notag
\end{align*}
The \em{rational weight} function $w_{l-a,a}$ is defined by the formula
\vvn.2>
\begin{align*}
w_{l-a,a}(t_1\lc t_l, z\); m_1, m_2)\,& {}=\!
\prod_{1\le u<v\le l}\,\frac{t_u\]-t_v}{t_u\]-t_v\]-1}
\;\prod_{u=1}^l\,\frac1{t_u\]+m_1}\;\x{}
\\[4pt]
& {}\>\x\,\Sym\biggl[\ \prod_{u=l-a+1}^l\)\frac{t_u}{t_u\]-z+m_2}\;
\prod_{1\le u<v\le l}\!\frac{t_u\]-t_v\]-1}{t_u\]-t_v}\ \biggr]
\\[-8pt]
\notag
\end{align*}
where \>$\Sym f(t_1\lc t_{l})\>=
\sum_{\sigma\in S_{l}}\,f(t_{\sigma_1}\lc t_{\sigma_{l}})$.
The \em{trigonometric weight} function $W_{l-b,b}$ is defined by the formula
\vvn.1>
\begin{align*}
W_{l-b,b}(t_1\lc t_l, z\); m_1, m_2)\,=\!
\prod_{1\le u<v\le l}\,\frac{\sink{t_u\]-t_v}}{\sink{t_u\]-t_v\]-1}}
\ \prod_{u=1}^l\,\frac{e^{\)-\pi\)i\)t_u/\ka}}{\sink{t_u\]+m_1}} &{}\;\x{}
\\[6pt]
{}\x\,\Sym\biggl[\ \prod_{u=l-b+1}^l\>
\frac{e^{\)\pi\)i\)z/\ka}\sin\)(\pi\)t_u/\ka)}{\sink{t_u\]-z+m_2}}\,
\prod_{1\le u<v\le l}\!\frac{\sink{t_u\]-t_v\]-1}}{\sink{t_u\]-t_v}}
&{}\ \biggr]\,.
\end{align*}
\vsk.4>
Integral \Ref{qmain} is defined for ${0<\Im\mu<2\pi}$. With respect to other
parameters we define integral \Ref{qmain} by analytic continuation from the
region where $m_1, m_2$ are complex numbers with negative real parts and
$\Re z=0$. In that case we put
\vvn.4>
$$
\dl_l(z\);m_1,m_2)\,=\,\lb\,(t_1\lc t_l)\in\C^{\>l}\ |\ \,
\Re t_u=\eps\,,\ \;u=1\lc l\,\rb
\vv.5>
$$
where $\eps$ is a positive number less than $\min\)(-\]\Re m_1\>, -\]\Re m_2)$.
In the considered region of parameters the integrand in \Ref{qmain} is well
defined on $\dl_{l_2}(z\); m_1,m_2)$ for any $a\),b$, the integral is
convergent and gives a meromorphic function of $z, m_1,\alb m_2$, see
\cite{TV1}\). It is also known that $I_{a,b}(z,\mu\);m_1,m_2,l_1,l_2)$
can be analytically continued to a nonnegative integer value of $m_2$,
if $a\),b$ are admissible with respect to $m_2, l_2$ at that point,
and the analytic continuation is given by the integral over
a suitable deformation of the imaginary plane
$\lb\,(t_1\lc t_{l_2})\in\C^{\>l_2}{\ |}\ \,\Re t_u=0$, $u=1\lc l_2\,\rb$,
see \cite{MuV}\).

\begin{remark}
There is an alternative way to describe the integrand of integral
\Ref{qmain}\), again writing it down as a product of three factors.
\vvgood
Namely, consider the functions
\begin{align*}
& \;\Xi_{\>l}(t_1\lc t_l, z,\mu\); m_1, m_2)\,=\,
(-\pi\ka)^{-\)l\)(\)l\)+3)/2}\,
\Exp{(\)\mu-\pi i\))\)\tsum_{u=1}^l\)t_u/\ka}\,\x{}
\\[4pt]
& {}\x\,\prod_{u=1}^l\,{\Gm(t_u/\ka)\,\Gm\bigl(\)-\)(m_1\]+t_u)/\ka\bigr)\,
\Gmk{t_u\]-z}\,\Gmk{z-m_2-t_u}}\,\x{}
\\[6pt]
& {}\x\!\prod_{1\le u<v\le l}\!(t_u\]-t_v)\,\sink{t_u\]-t_v}\,
\Gmk{t_u\]-t_v\]+1}\,\Gmk{t_v\]-t_u\]+1}\,,
\\
\gooddb\nn3>
\end{align*}
\vsk-2>
$$
p_{l-a,a}(t_1\lc t_l, z\);m_2)\,={}
\,\Sym\biggl[\ \prod_{u=1}^{l-a}\,(t_u\]-z+m_2)\!\prod_{u=l-a+1}^l\!\!t_u
\prod_{1\le u<v\le l}\!\frac{t_u\]-t_v\]-1}{t_u\]-t_v}\ \biggr]\,,
\ifUS\else\kern-.5em\fi
\vv.4>
$$
\begin{align*}
& P_{l-b,b}(t_1\lc t_l, z\);m_2)\,=\,
\Exp{\pi i\)bz/\ka}\,\x{}
\\[4pt]
&\!\}{}\x\,\Sym\biggl[\ \prod_{u=1}^{l-b}\,\sink{t_u\]-z+m_2}\!
\prod_{u=l-b+1}^l\!\sin\)(\pi\)t_u/\ka)
\prod_{1\le u<v\le l}\!\frac{\sink{t_u\]-t_v\]-1}}{\sink{t_u\]-t_v}}
\ \biggr]\,.
\\[-12pt]
\end{align*}
Then
\vvn-.6>
\begin{align*}
\Phi_{l_2}(t,z,\mu\); m_1, m_2)\,
w_{l_2-a,a}(t, z\);m_1, m_2)\,W_{l_2-b,b}(t, z\);m_1, m_2)\,={}\!\] &
\\[8pt]
{}=\,\Xi_{l_2}(t,z,\mu\); m_1, m_2)\,
p_{l_2-a,a}(t, z\); m_2)\,P_{l_2-b,b}(t, z\); m_2)\,. &
\\[-8pt]
\end{align*}
\par
In this description of the integrand the function $\Xi_l$ is such that
it contains all the poles of the integrand and has no zeros anywhere but on
the shifted diagonals: $(t_u\]-t_v)\in\ka\>\Z$, the functions $p_{l-a,a}$
are polynomials in $t_1\lc t_l, z, m_2$, and the functions $P_{l-b,b}$
are trigonometric polynomials in $t_1\lc t_l, z, m_2$.
\end{remark}

\subsection{Hypergeometric integrals}

Let $\ka$ be a positive number.
Let $m_1, l_1$ be complex numbers and $m_2, l_2$ nonnegative integers
such that
\vvn.2>
$$
m_1+\>m_2\,=\,l_1+\>l_2\,.
\vv.2>
$$
\par
For a pair of admissible numbers $a, b$ we define a function
$J_{a,b}(z,\mu\); m_1, m_2, l_1, l_2)$ of complex variables $z,\mu$.
The function is defined as an \$l_2\)$-dimensional hypergeometric integral:
\vvn.4>
\begin{equation}
\label{main}
J_{a,b}(z,\mu\); m_1, m_2, l_1, l_2)\,=\!\int_{\gm_{\)l_2-b,\)b}(z)\!}\!\!\!
\Psi_{l_2}(t,z,\mu\);m_1,m_2)\,g_{\)l_2-a,\)a}(t, z)\,dt^{\>l_2}\,.
\end{equation}
\vsk.2>\nt
Here ${t=(t_1\lc t_{l_2})}$, \,${dt^{\>l_2} = dt_1\ldots\)dt_{l_2}}$.
The functions $\Psi_{l_2}(t,z,\mu\);m_1, m_2)$ and
$g_{l_2-a,a}(t,\alb z)$ are defined below. The \$l_2$-dimensional
integration contour $\gm_{l_2-b,\)b}(z)$ lies in $\C^{\>l_2}$ and
is also defined below.

The \em{master} function $\Psi_l$ is defined by the formula
\vvn.4>
\begin{align*}
\Psi_l &{} (t_1\lc t_l, z,\mu\); m_1, m_2)\,={}
\\[6pt]
&{}\]=\,\prod_{u=1}^l\,t_u^{\)(\)\mu\)+\)m_1+\)m_2-2\)l\)+\)1)/\ka}\>
(1-t_u)^{-m_1/\ka}\>(z-t_u)^{-m_2/\ka}\]
\prod_{1\le u<v\le l}\}(t_u\]-t_v)^{2/\ka}\,.\kern-1em
\\[-9pt]
\gooddb
\end{align*}
The \em{weight} function $g_{l-a,\)a}$ is defined by the formula
\vvn.3>
$$
g_{\)l-a,\)a}(t_1\lc t_l, z) \,=\,
\Sym \biggl[\ \prod_{u=1}^{l-a}\,\frac 1{1-t_u}\
\prod_{u=l-a+1}^l\,\frac 1{z-t_u}\ \biggr]\,.
$$
\vsk.3>
We define the integral in \Ref{main} by analytic continuation from the region
\vvn.3>
\begin{equation}
\label{region}
z\ne 0\,,\qquad 0<\arg\)z<2\)\pi\,,\qquad \Re\mu\ll 0\,.
\end{equation}
\vsk.4>\nt
In that region the integration contour
$\gm_{\)l-b,\)b}(z)$ is shown in the picture.

\vsk>
\vbox{\begin{center}
\begin{picture}(185,135)
\put(45,15){\oval(10,10)[l]}
\put(45,15){\circle*{3}}
\put(45,10){\vector(1,0){152}}
\put(45,20){\line(1,0){152}}
\qbezier[16](51,5)(120,5)(195,5)
\qbezier[16](51,25)(120,25)(195,25)
\qbezier[4](35,15)(35,25)(45,25)
\qbezier[4](45,5)(35,5)(35,15)
\put(45,15){\oval(30,30)[l]}
\put(45,0){\vector(1,0){152}}
\put(45,30){\line(1,0){152}}

\put(10,45){\circle*{3}}
\qbezier(6,41)(2,45)(5.8,49)
\qbezier(6,41)(10,38.1)(14,41.2)
\put(6,49){\line(1,1){69}}
\put(14,41){\vector(1,1){69}}
\qbezier(-1,34)(-12,45)(-1.2,56)
\qbezier(-1,34)(10,24.4)(21,34.3)
\put(-1,56){\line(1,1){69}}
\put(21,34){\vector(1,1){69}}
\qbezier[10](21,41)(53,73)(85,105)
\qbezier[10](6,56)(38,88)(70,120)
\qbezier[4](17,37)(10,31)(2,37)
\qbezier[4](2,52)(-5,45)(2,37)

\put(17,13){$1$}
\put(202,30){${t_l}$}
\put(202,20){${t_{b+1}}$}

\put(-19,43){$z$}
\put(69,130){${t_b}$}
\put(78,121){${t_1}$}

\put(-35,13){$0$}
\put(-20,15){\circle*{3}}
\end{picture}
\vsk1.3>
Picture~1. Integration contour $\gm_{\)l-b,b}(z)$.
\end{center}}

\vsk.6>
\nt
It has the form ${\gm_{\)l-b,\)b}(z)\,=\,\lb\,(t_1\lc t_l) \in \C^{\>l}
\ |\ \,t_u\in\Cc_u\,,\ \;u = 1\lc l\;\rb}$.
Here $\Cc_u$, $u=1\lc l$, are non-intersecting oriented loops in $\C$.
The first $b$ loops start at infinity in the direction of $z$, go around
$z$, and return to infinity in the same direction. For $1\le u<v\le b$
the loop $\Cc_u$ lies inside the loop $\Cc_v$. The last $l-b$ loops start at
infinity in the real positive direction, go around $1$, and return to infinity
in the same direction. For $b+1\le u<v\le l$ the loop $\Cc_u$ lies inside
the loop $\Cc_v$.
\vsk.1>
If $z,\mu$ are in region \Ref{region}, we fix a univalued branch of the master
function $\Psi_l$ over $\gm_{\)l-b,\)b}(z)$ by fixing the arguments of all its
factors. Namely, we assume that at the point of $\gm_{\)l-b,\)b}(z)$ where all
numbers $t_1/\]z\>\lc\)t_b/\]z\>,\)t_{b+1}\>\lc\)t_l$ belong to
\vvn.5>
$(0\>,1)$ we have
$$
\Line{\arg\>t_u\in[\)0\>,2\)\pi)\,,\hfil \arg\)(1-t_u)\in(-\)\pi\),\pi)\,,\hfil
\arg\)(z-t_u)\in(0\>,2\)\pi)\,,\hfil \arg\)(t_u\]-t_v)\in[\)0\>,2\)\pi)\,,}
\vv.5>
$$
for $u=1\lc l$, $v=u+1\lc l$.
Integral \Ref{main} is convergent in region \Ref{region}.

\subsection{Main result}

\begin{theorem}\label{first}
Let $\ka$ be a generic positive number. Let $m_1, l_1$ be complex numbers
and $m_2, l_2$ nonnegative integers, such that $m_1+m_2 = l_1+l_2$.
Let \>${0<\Im\mu<2\pi}$. Then for any $a,b=1\lc\min\)(m_2\),l_2)$ we have
\ifUS\vvn.4>\else\vvn.5>\fi
\begin{align}
\label{I=J}
C_b(m_1,m_2,l_1,l_2)\;{}& I_{a,b}(z,\mu\);m_1,m_2,l_1,l_2)\,=\!\]{}
\\[7pt]
\ifUS\else\nn2>\fi
{}=\,D_b(m_1, m_2,l_1,l_2)\,E_b(l_1,l_2)\,{}& X(z\);m_1,m_2)\,
Y(\mu\);m_1,m_2,l_1,l_2)\,J_{a,b}(\emu,z\);l_1,l_2,m_1,m_2)
\notag
\\[-14pt]
\gooddb\ifUS\else\nn1>\fi
\notag
\end{align}
where
\ifUS\vvn-.4>\else\vvn-.3>\fi
\begin{align*}
C_b(m_1,m_2,l_1,l_2)\,=\;\frac{(2\)\pi i)^{-\)l_2}}{l_2!\>(l_2-b)!\,b!}\,
\prod_{j=0}^{l_2-b-1}\sink{m_1-j}\,\prod_{j=0}^{b-1}\,\sink{m_2-j}\,\x{}
\\[5pt]
{}\x\,\prod_{j=0}^{l_2-1}\,\frac{\Gm(1+1/k)\>\Gm\bigl(1+(m_1-j)/\ka\bigr)}
{\Gm\bigl(1+(j+1)/\ka\bigr)}\; & ,
\end{align*}
\vvn.2>
\begin{align*}
D_b(m_1,m_2,l_1,l_2)\,=\,(2\)i)^{-\)m_2}\,
\prod_{j=0}^{m_2-b-1}\!\frac1{\sink{j+1}}\,\;
\prod_{j=0}^{b-1}\,\frac1{\sink{j+1}}\;\x{}
\\[5pt]
{}\x\,\prod_{j=0}^{m_2-1}\;\frac{\Gm\bigl(1+(l_1-j)/\ka\bigr)}
{\Gm(-1/\ka)\>\Gm\bigl(1+(j+1)/\ka\bigr)}\; & ,
\end{align*}
\vvn.4>
$$
E_b(l_1,l_2)\,=\,
\expk{b^{\)2}-\)(b\>-\)l_2)\>(l_1+\)l_2)-\)l_2\)(l_2-1)/2}\,,
$$
\vvn0>
\begin{equation}
\label{Xt}
X(t\);m_1,m_2)\,=\>\prod_{j=0}^{m_2\)-1}\,\frac{\Gmk{j-m_1-t}}{\Gmk{j+1-t)}}\;,
\end{equation}
\vvn.2>
\begin{equation}
\label{Y}
Y(\mu\);m_1,m_2,l_1,l_2)\,=\,
e^{\>\mu\>l_2(l_2\)-\)2\)m_1\)-1)/2\ka}\>(1-\emu)^{(l_1+1)\)l_2/\ka}\,,
\quad \arg\)(1-\emu)\in(\)-\)\pi\),\pi)\,.\ifUS\else\kern-1em\fi
\end{equation}
\vsk.3>
\end{theorem}
\begin{remark}
There is a similar theorem establishing an equality of suitable
hypergeometric integrals of different dimensions, see \cite{TV4}\).
The factor $D_b(m_1,m_2,l_1,l_2)$ in the present paper corresponds to
the normalization factor $C_b(l_1,l_2,m_1,m_2)$ in \cite{TV4} and contains
the same product of sines and gamma-functions.
\end{remark}
\begin{example}
Let $m_2=l_2=1$. In this case formula \Ref{I=J} becomes the classical equality
of two integral representations of the Gauss hypergeometric function $_2F_1$.
For instance, if $a=b=0$, then taking ${\al=-\)m_1/\ka}$,
${\bt=-\)(z+m_1)/\ka}$, ${\gm=(1-z-m_1)/\ka}$,
after simple transformations one gets:
\vvn.6>
\begin{align*}
& \frac1{2\pi i}\,\int_{-\)i\)\infty\)-\eps}^{+\)i\)\infty\)-\)\eps}\!
e^{\>s\)(\mu-\pi i)}\,\Gm(-\)s)\>\Gm(s+\al)\;\frac{\Gm(s+\bt)}{\Gm(s+\gm)}\ ds
\,={}
\\[2pt]
&\ \ {}=\,(1-\emu)^{\gm-\al-\bt}\;\frac{\Gm(\bt)}{\Gm(\gm-\al)}\,\,
\int_{1\,}^{+\infty} t^\bt\>(t-1)^{\al-1}\>(t-\emu)^{\bt-\gm}\>dt\,={}
\\[2pt]
&\ \ \quad{}=\;\frac{\Gm(\bt)}{\Gm(\gm-\al)}
\,\,\int_{\!0}^{1\!}u^{\al-1}\>(1-u)^{\gm-\al-1}\>(1-u\)\emu)^{-\)\bt}\>du\,=
\;\frac{\Gm(\al)\>\Gm(\bt)}{\Gm(\gm)}\ {}_2F_1(\al\),\bt\);\gm\);\emu)\,.
\\[-10pt]
\gooddb
\end{align*}
Here it is assumed that ${\Re\gm>\Re\al>0}$, ${\Re\bt>0}$ and
${0<\eps<\min\)(\)\Re\al\>,\Re\bt\))}$. The second equality is obtained
by the change of integration variable $u\>=\>(t-1)/(t-\emu)$.
\end{example}
Theorem~\ref{first} claims that an \$l_2$-dimensional \q-hypergeometric
integral equals an \$m_2$-dimen\-sional hypergeometric integral up to an
explicit factor. Note that in the first integral the numbers $m_1, m_2$ are
shifts of arguments of the gamma-functions entering its \q-master function,
while in the second integral $m_2$ is its dimension and $m_1$ is not present
explicitly.
\vsk.1>
It is well known that studying asymptotics of integrals with respect to their
dimension is an interesting problem appearing, for instance, in the theory of
orthogonal polynomials and in matrix models. The duality of the theorem allows
us to study asymptotics of integrals with respect to their dimension. Namely,
assume that in the 4\)-tuple $(m_1,m_2,l_1,l_2)$ the nonnegative integer $l_2$
tends to infinity while the numbers $m_1$ and $m_2$ remain fixed.
Then $I_{a,b}(z,\mu\); m_1, m_2, l_1, l_2)$ is a \q-hypergeometric integral
of growing dimension $l_2$, whose master function has fixed shifts $m_1,m_2$.
At the same time, $J_{a,b}(\emu,z\);l_1,l_2,m_1,m_2)$ is a hypergeometric
integral of the fixed dimension $m_2$, whose master function has growing
exponents $l_1,l_2$. The asymptotics of $J_{a,b}(\emu,z\);l_1,l_2,m_1,m_2)$ can
be calculated using the steepest descent method. An example of such calculation
is given in \cite{TV4}\).
\vsk.1>
Similarly one can assume that in the 4\)-tuple $(m_1, m_2, l_1, l_2)$
the nonnegative integer $m_2$ tends to infinity while the numbers $l_1$ and
$l_2$ are fixed. Then $J_{a,b}(\emu,z\);l_1,l_2,\alb m_1,m_2)$ is
a hypergeometric integral of growing dimension $m_2$, whose master function
has fixed exponents $l_1, l_2$. On the other hand,
$I_{a,b}(z,\mu\);m_1,m_2,l_1,l_2)$ is a \q-hypergeometric integral of
the fixed dimension $l_2$, whose master function has growing shifts $m_1, m_2$.
The asymptotics of $I_{a,b}(z,\mu\);m_1,m_2,l_1,l_2)$ in principle can be
calculated using the Stirling formula for asymptotics of the gamma-function.
\vsk.2>
To prove Theorem~\ref{first} we show that the matrices:
\vvn.4>
\begin{equation}
\label{Imat}
\bigl(I_{a,b}(z,\mu\);m_1,m_2,l_1,l_2)\bigr)_{0\le a,\)b\)\le\min\)(m_2,\)l_2)}
\end{equation}
\vsk-.2>\nt
and
\vvn-.3>
$$
X(z\);m_1,m_2)\>\bigl(J_{a,b}(\emu\},z\);l_1,l_2,m_1,m_2)\bigr)
_{0\le a,\)b\)\le\min\)(m_2,\)l_2)}
\vv.4>
$$
satisfy the same system of first order linear difference equations with
respect to $z$, see Corollary~\ref{Xz} and Theorems~\ref{I-qKZ}, \ref{J-KZ}.
Studying asymptotics of those matrices as $\Re z\to-\)\infty$ allows
us to compute the connection matrix, thus proving Theorem~\ref{first}.
\vsk.2>
The fact that both matrices satisfy the same system of difference equations
is based on the duality of the \qKZ/ and dynamical equations for $\gl_k$ and
$\gl_n$ \cite{TV3}\), which is a generalization of the duality between the
rational differential \KZ/ and dynamical equations observed in \cite{TL}\).
Namely, in \cite{TV2} a system of dynamical difference equations was
introduced. It is proved there that the dynamical difference equations are
compatible with the trigonometric \KZ/ differential equations. In \cite{MV}
hypergeometric solutions of the trigonometric \KZ/ and dynamical difference
equations were presented. On the other hand, \q-hypergeometric solutions of the
rational \qKZ/ difference equations equations were constructed in \cite{TV1}.
It was shown in \cite{TV3} that the system of dynamical difference equations
for $\gl_k$ and the system of rational \qKZ/ equations for $\gl_n$ are
naturally transformed into each other under the $(\gl_k\>,\gl_n)$ duality.
In this way one gets two sets of solutions of the same system of equations,
and in principle, these two sets of solutions can be identified.
Theorem~\ref{first} is a realization of this idea for the case of $k=n=2$.
We will discuss the case of an arbitrary pair $k, n$ in a separate paper.
\vsk.1>
Let us make an additional remark. In \cite{TV3} it was introduced a system
of dynamical differential equations which is transformed under the
$(\gl_k\>,\gl_n)$ duality to the trigonometric \KZ/ differential equations.
It is proved in \cite{TV5} that the dynamical differential equations are
compatible with the rational \qKZ/ difference equations. It is shown in
\cite{TV6} that the \q-hypergeometric solutions of the \qKZ/ equations satisfy
the dynamical differential equations. In the present case of $k=n=2$ this means
that the matrices \Ref{Imat} and
\vvn.4>
\begin{equation}
\label{YJ}
Y(\mu\);m_1,m_2,l_1,l_2)\>\bigl(J_{a,b}(\emu,z\);l_1,l_2,m_1,m_2)\bigr)
_{0\le a,\)b\)\le\min\)(m_2,\)l_2)}
\end{equation}
\vsk.2>\nt
satisfy the same system of first order linear differential equations
with respect to $\mu$. Therefore, the connection matrix of \Ref{Imat} and
\Ref{YJ} does not depend on $\mu$, which agrees with formula \Ref{I=J}\).
Other factors in \Ref{I=J} also have natural explanation in terms of
the \q-deformation of the $(\gl_k\>,\gl_n)$ duality.
\vsk.2>
The rest of the paper is as follows. In Section~2 we discuss
the $(\glt\>,\glt)$ duality for the \qKZ/ and dynamical equations. In Section~3
we describe their hypergeometric solutions and prove Theorem~\ref{first}.
In Section~4 we give some facts about the dynamical differential equations.

\vsk.2>
The authors thank Y\]\&Markov for useful discussions.

\section{The $(\glt\>,\glt)$ duality for \KZ/ and dynamical equations}
In this section we follow the exposition in \cite{TV3}.

\subsection{The trigonometric \KZ/ and associated dynamical difference
equations}

Let $E_{ij}$, $i,j=1,2$, be the standard generators of the Lie algebra $\glt$.
Let $\hg=\C\>{\cdot}E_{11}\oplus\)\C\>{\cdot}E_{22}$ be the Cartan subalgebra,
while $E_{12}$ and $E_{21}$ are respectively the positive and negative root
vectors.
\vsk.2>
The trigonometric \r-matrix $r(z)\in\glt^{\)\ox2}$ is defined by the formula
\vvn.2>
$$
r(z)\,=\,\frac1{z-1}\,\bigl((z+1)\>(E_{11}\ox\]E_{11}+E_{22}\ox\]E_{22})/2\>+
\)z\)E_{12}\ox\]E_{21}\)+\)E_{21}\ox\]E_{12}\bigr)
$$
\vsk.2>
Fix a non-zero complex number $\ka$. The \em{trigonometric \KZv/} (\KZ/)
\em{operators} $\nabla_1\)\lc\]\nabla_n$ are the following differential
operators with coefficients in $U(\glt)^{\ox\)n}\}$ acting on functions
of complex variables $z_1\lc z_n,\>\la_1,\la_2$:
\vvn.2>
\begin{equation}
\label{KZo}
\nabla_a(z_1\lc z_n,\la_1,\la_2)\,=\,
\ka\)z_a\>\frac\partial{\partial z_a\]}\,-\)
\sum_{i=1}^2\,\bigl(\la_i\>-\>\frac12\>\tsum_{b=1}^n E_{ii}^{\)(b)}\)\bigr)
\>E_{ii}^{(a)}\)-\)\sum_{\fratop{b=1}{b\ne a}}^n\>r(z_a/z_b)\vpb{(a\)b)}\).
\end{equation}
Here
$E_{ii}^{\)(b)}=\)1\lox{}\%{E_{ii}}_{\Clap{\sss\text{\=$b$-th}}}\]{}\lox 1$,
\>and the meaning of $r(z_a/z_b)\vpb{(a\)b)}\}$ is similar. It is known that
the operators $\nabla_1\)\lc\]\nabla_n$ pairwise commute.
\vsk.2>
Let $V_1\lc V_n$ be \$\glt$-modules. The \em{trigonometric \KZ/ equations} for
\vvn.06>
a function $U(z_1\lc z_n,\la_1,\la_2)$ with values in $V_1\lox V_n$ are
\vvn.4>
\begin{equation}
\label{KZ}
\nabla_a(z_1\lc z_n,\la_1,\la_2)\,U(z_1\lc z_n,\la_1,\la_2)\,=\,0\,,
\qquad a=1\lc n\,.\kern-1em
\end{equation}
\goodbreak
\vsk.4>
Let $T_u$ be the difference operator acting on functions $f(u)$ by the formula
\vvn.3>
$$
(T_u f)(u)\,=\,f(u+\ka)\,.
$$
\vsk.2>
Introduce a series $B(t)$ depending on a complex variable $t$:
\vvn.4>
$$
B(t)\,=\,1\>+\>\sum_{s=1}^\infty\,E_{21}^s\)E_{12}^s\,
\prod_{j=1}^s\,\frac1{j\>(t-E_{11}+E_{22}-j\))}\;.
\vv.4>
$$
For any \$\glt$-module $V\]$ with a locally nilpotent action of $E_{12}$
and finite-dimensional weight subspaces the series $B(t)$ has a well-defined
action in any weight subspace $V[\)\mu\)]\sub V\]$ as a rational
\$\End\bigl(V[\)\mu\)]\bigr)\)$-valued function of $t$.
\vsk.2>
Let $V_1\lc V_n$ be \$\glt$-modules as above. Introduce the \em{dynamical
difference operators} $Q_1\>,\)Q_2$ acting on \$V_1\lox V_n\)$-valued functions
of complex variables $z_1\lc z_n,\>\la_1,\la_2$ by the formulae
\vvn-.2>
\begin{gather}
Q_1(z_1\lc z_n,\la_1,\la_2)\,=\,\bigl(B(\la_1-\la_2)\bigr)^{-1}\,
\prod_{a=1}^n\>z_a^{\)-E_{11}^{(a)}}\,T_{\la_1}\,,
\\
Q_2(z_1\lc z_n,\la_1,\la_2)\,=\,
\prod_{a=1}^n\>z_a^{\)-E_{22}^{(a)}}B(\la_1-\la_2-\ka)\,T_{\la_2}\,.
\notag
\\[-13pt]
\notag
\end{gather}
One can see that the operators $Q_1\>,\)Q_2$ commute.
\begin{proposition}[\)\cite{TV2}\)]
One has $[\)\nabla_a\>,Q_i\)]\)=\)0$ for all $a=1\lc n$ and $i=1,2$.
\end{proposition}
The \em{dynamical difference equations} associated with the trigonometric
\KZ/ equations for a function $U(z_1\lc z_n,\la_1,\la_2)$ with values in
$V_1\lox V_n$ are
\vvn.4>
\begin{equation}
\label{q-dyn}
Q_i(z_1\lc z_n,\la_1,\la_2)\,U(z_1\lc z_n,\la_1,\la_2)\,=\,
U(z_1\lc z_n,\la_1,\la_2)\,,\qquad i=1,2\,.\kern-1em
\end{equation}
\vsk.4>
The trigonometric \KZ/ and dynamical difference operators preserve the weight
decomposition of $V_1\lox V_n$. Thus the \KZ/ and dynamical equations can be
considered as equations for a function $U(z_1\lc z_n,\la_1,\la_2)$ taking
values in a given weight subspace of $V_1\lox V_n$.

\subsection{The \qKZ/ equations}
Let $V, W\}$ be irreducible highest weight \$\glt\)$-modules with highest
weight vectors $v\), w$, respectively. There is a unique rational function
$R_{VW}(t)$ taking values in $\End(V\]\ox W)$ such that
\vvn.6>
$$
\bigr[\)R_{VW}(t)\>,\>g\otimes 1+1\otimes g\>\bigr]\,=\,0\qquad
\text{for any}\quad g\in\glt\,,
\vv.7>
$$
$$
R_{VW}(t)\>
\bigl(E_{21}\ox E_{11}\)+E_{22}\ox E_{21}\)+\)t\>E_{21}\ox 1\)\bigr)\,={}
\,\bigl(E_{11}\ox E_{21}\)+E_{21}\ox E_{22}\)+\)t\>E_{21}\ox 1\)\bigr)
\>R_{VW}(t)\,,
\vv.8>
$$
$$
R_{VW}(t)\,v\ox w\,=\,v\ox w\,.
\vv.5>
$$
The function $R_{VW}(t)$ is called the rational \R-matrix for the tensor
product $V\]\ox W\}$. It comes from the representation theory of the Yangian
$Y(\glt)$.
\vsk.2>
Let $V_1\lc V_n$ be irreducible highest weight \$\glt$-modules.
Introduce the \em{\qKZ/ difference operators} $Z_1\lc Z_n$ acting
on \$V_1\lox V_n\)$-valued functions of complex variables
$z_1\lc z_n,\>\la_1,\la_2$ by the formula
\vvn.3>
\begin{align}
\label{qKZo}
Z_a(z_1\lc z_n, {}&\la_1,\la_2)\,=\,\bigl(\)R_{an}(z_a\]-z_n)\ldots
R_{a,\)a+1}(z_a\]-z_{a+1})\)\bigr)\vpb{-1}\x{}
\\[5pt]
{}\x\, {}& \la_1^{\)-E_{11}^{(a)}}\)\la_2^{\)-E_{22}^{(a)}}\)
R_{1a}(z_1-z_a\]-\ka)\ldots R_{a-1,\)a}(z_{a-1}-z_a\]-\ka)\,T_{z_a}\,.
\notag
\\[-12pt]
\notag
\end{align}
These operators are called the \em{\qKZ/ operators}. It is known that they
pairwise commute \cite{FR}\). The difference equations
\vvn.5>
\begin{equation}
\label{qKZ}
Z_a(z_1\lc z_n,\la_1,\la_2)\,U(z_1\lc z_n,\la_1,\la_2)\,=\,
U(z_1\lc z_n,\la_1,\la_2)\,,\quad\ \ a=1\lc n\,,\kern-1em
\end{equation}
\vsk.3>\nt
for a \$V_1\lox V_n\)$-valued function $U(z_1\lc z_n,\la_1,\la_2)$ are called
the \em{\qKZ/ equations}.
\vsk.2>
The \qKZ/ operators preserve the weight decomposition of $V_1\lox V_n$.
Thus the \qKZ/ and dynamical differential equations can be considered
as equations for a function $U(z_1\lc z_n,\la_1,\la_2)$ taking values
in a given weight subspace of $V_1\lox V_n$.

\subsection{The duality}
For a complex number $m$, denote $M_m$ the Verma module over $\glt$
with highest weight $(m,0)$ and highest weight vector $v_m$.
The vectors $E_{2,1}^d\)v_m$, $d\in\Zp$, form a basis in $M_m$.
\vsk.1>
For a non-negative integer $m$, denote $L_m$ the irreducible $\glt$-module
with highest weight $(m,0)$ and highest weight vector $v_m$. The vectors
$E_{2,1}^d\)v_m$, $d=0\lc m$, form a basis in $L_m$.
\vsk-.4>
\goodbreak
\vsk.6>
Let $m_1,l_1$ be complex numbers and $m_2, l_2$ non-negative integers
such that ${m_1+m_2}={l_1+l_2}$. Consider the weight subspace
${(M_{m_1}\}\ox L_{m_2})}[\)l_1\),l_2]$ of the tensor product
${M_{m_1}\}\ox L_{m_2}}$. The weight subspace has a basis
\vvn.3>
\begin{equation}
\label{Fbasis}
F^a(m_1,m_2,l_1,l_2)\,=\,
\frac1{(l_2-a)!\,a\)!}\;E_{2,1}^{l_2-a}\)v_{m_1}\ox E_{2,1}^a\)v_{m_2}\,,
\qquad a=0\lc\min\)(m_2, l_2)\,.\!\!
\end{equation}
\vsk.2>\nt
There is a linear isomorphism
\vvn.3>
\begin{align}
\label{phiso}
\phi\):\)(M_{m_1}\ox L_{m_2})[\)l_1\),l_2]\, &{}\to\,
(M_{l_1} \ox L_{l_2})[m_1,m_2]\,,
\\[7pt]
F^a(m_1,m_2,l_1,l_2)\, &{}\mapsto\> F^a(l_1,l_2,m_1,m_2)\,.
\notag
\end{align}
\vsk.5>
\begin{theorem}[\)\cite{TV3}\)]
The isomorphism $\phi$ transforms the \qKZ/ operators acting in\ifUS\else\nl\fi
$(M_{m_1}\}\ox L_{m_2})[\)l_1\),l_2]$ into the dynamical difference
operators acting in $(M_{l_1}\}\ox L_{l_2})[m_1,m_2]$. More precisely, we have
\vvn.4>
\begin{align*}
& \phi\;Z_1(z_1,z_2,\la_1,\la_2)\,=\,
\bigl(G(z_1-z_2;m_1,m_2)\bigr)\vpb{-1}\)Q_1(\la_1,\la_2,z_1,z_2)\;\phi\,,
\\[6pt]
& \phi\;Z_2(z_1,z_2,\la_1,\la_2)\,=\,
G(z_1-z_2-\ka;m_1,m_2)\,Q_2(\la_1,\la_2,z_1,z_2)\;\phi\,,
\\[-12pt]
\end{align*}
where
\vvn-.4>
$$
G(t\);m_1,m_2)\,=\,\prod_{j=0}^{m_2-1}\>\frac{t+j-m_1}{t+j+1}\;.
$$
\end{theorem}
Let $S(t)$ be any solution of the equation
\vvn.3>
$$
S(t+\ka)\,=\,G(t\);m_1,m_2)\,S(t)\,.
\vv.3>
$$
For instance, one of solutions is given by $S(t)=X(-\)t\);m_1,m_2)$,
cf.~\Ref{Xt}.
\begin{corollary}
\label{Xz}
Let an \$(M_{l_1}\]\ox L_{l_2})[m_1,m_2]$-valued function
$U(z_1,z_2,\la_1,\la_2)$ solve the dynamical difference equations\/{\rm:}
\vvn.4>
$$
Q_a(z_1,z_2,\la_1,\la_2)\,U(z_1,z_2,\la_1,\la_2)\,=\,U(z_1,z_2,\la_1,\la_2)\,,
\qquad a=1,2\,.
\vv.4>
$$
Then the \$(M_{m_1}\]\ox L_{m_2})[\)l_1\),l_2]$-valued function
\vvn.4>
$$
\Ut(z_1,z_2,\la_1,\la_2)\,=\,
S(z_1-z_2)\,\phi^{-1}\bigl(U(\la_1,\la_2,z_1,z_2)\bigr)
\vv.4>
$$
solves the \qKZ/ equations\/{\rm:}
\vvn.4>
$$
Z_a(z_1,z_2,\la_1,\la_2)\,\Ut(z_1,z_2,\la_1,\la_2)\,=\,
\Ut(z_1, z_2,\la_1,\la_2)\,, \qquad a=1,2\,.
$$
\end{corollary}
\vsk.8>
More facts on the $(\glt\>,\glt)$ duality for \KZ/ and dynamical equations
are given in Section~4.

\section{Hypergeometric solutions}
\subsection{Hypergeometric solutions of the \qKZ/ equations}
For any $b=0\lc\min\)(m_2,l_2)$ define
an \$(M_{m_1}\]\ox L_{m_2})[\)l_1\),l_2]$-valued function
\vvn.4>
$$
\Ib_b(z,\mu\);m_1,m_2,l_1,l_2)\,=\>
\sum_{a=0}^{l_2}\,I_{a,b}(z,\mu\); m_1,m_2,l_1,l_2)
\;\frac1{(l_2-a)!\,a\)!}\,E^{l_2-a}_{2,1}\)v_{m_1} \ox E^{a}_{2,1}\)v_{m_2}\,.
\vv.3>
$$
Here the actual range of summation is until $a=\min\)(m_2,l_2)$,
since $E^{a}_{2,1}\)v_{m_2}=0$ for $a>m_2$.

\begin{theorem}
\label{I-qKZ}
Let \>$\la_i=e^{\>\mu_i}\}$, $i=1,2$.
For any $b=0\lc\min\)(m_2,l_2)$ the function
\vvn.5>
\begin{align}
\label{Ub}
U_b(z_1,z_2,\la_1,\la_2)\, &{}=\,
e^{\)(\)\mu_1(m_1z_1+\>m_2\)z_2\)-\>(m_1^2+\>m_2^2\))/2)
\>-\>(\mu_1\]-\)\mu_2)(l_2z_1\]+\>l_2/2))/\ka}\,\x{}
\\[7pt]
&\>{}\x\>(1-e^{\>\mu_2-\mu_1})^{-\)l_2/\ka}\>
\Ib_b(z_2-z_1,\mu_2\)-\mu_1;m_1,m_2,l_1,l_2)\kern-1em
\notag
\\[-12pt]
\notag
\end{align}
is a solution of the \qKZ/ equations \Ref{qKZ} with values in
$(M_{m_1}\]\ox L_{m_2})[\)l_1\),l_2]$. Moreover, if $\la_1/\la_2$ is not real,
than any solution of that \qKZ/ equations is a linear combination of functions
$U_b(z_1,z_2,\la_1,\la_2)$ with coefficients being \$\ka$-periodic functions
of $z_1,z_2$.
\end{theorem}
\nt
The theorem is a direct corollary of the construction of \q-hypergeometric
solutions of the \qKZ/ equations given in \cite{TV1}\), \cite{MuV}\).
\vsk.5>
We describe asymptotics of the integral $I_{a,b}(z,\mu\);m_1,m_2,l_1,l_2)$
as $\Re z\to-\)\infty$ and $\mu$ is fixed. For a positive number $\ka$,
complex numbers $m\>,\mu$, ${\Im\mu\in(0\>,2\)\pi)}$, and a nonnegative
integer $l$ consider the Selberg-type integral
\vvn.2>
\begin{align}
\label{q-selberg}
A_l(\mu\);m)\,=\int_{\dl_l(m)}\!\Exp{(\)\mu-\pi i\))\tsum_{u=1}^l s_u}\,
\prod_{u=1}^l\,\Gm(s_u)\,\Gm(\)-\)s_u\]-m/\ka)\,\x{}\} &
\\[-2pt]
{}\x\,\prod_{\fratop{u,v=1}{u\ne v}}^l
\frac{\Gm(s_u\]-s_v\]+1/\ka)}{\Gm(s_u\]-s_v)}\ ds^{\)l}\, &.
\notag
\\[-14pt]
\notag
\end{align}
The integral is defined by analytic continuation from the region where $\Re m$
is negative. In that case
\vvn.1>
$$
\dl_l(m)\,=\,\lb\,(s_1\lc s_l)\in\C^{\>l}\ |
\ \,\Re s_u=\)\Re m/2\,,\ \;u=1\lc l\,\rb\,.
\vv.5>
$$
In the considered region of parameters the integrand in \Ref{q-selberg} is well
defined on $\dl_l(m)$ and the integral is convergent, see \cite{TV1}\).
The formula for $A_l(m)$ is well known,
\vvn.2>
\ifUS
$$
A_l(\mu\);m)\,=\,(2\)\pi i)^l\,e^{\)(\mu\)-\pi i)\)(l-1-2\)m)\>l/2\ka}\,
(1-\emu)^{\)l\)(m-l+1)/\ka}\,\prod_{j=0}^{l-1}\,
\frac{\Gm\bigl(1+(j+1)/\ka)}{\Gm(1+1/\ka)}\;\Gmk{j-m}\,,
\vv.5>
$$
\else
\begin{align*}
& \,A_l(\mu\);m)\,={}
\\
& {}=\,(2\)\pi i)^l\,e^{\)(\mu\)-\pi i)\)(l-1-2\)m)\>l/2\ka}\,
(1-\emu)^{\)l\)(m-l+1)/\ka}\,\prod_{j=0}^{l-1}\,
\frac{\Gm\bigl(1+(j+1)/\ka)}{\Gm(1+1/\ka)}\;\Gmk{j-m}\,,
\\[-12pt]
\end{align*}
\fi
where $\arg\)(1-\emu)\in(-\)\pi\),\pi)$, see, for example, \cite{TV1}\).
\begin{remark}
Other versions of Selberg-type integrals see in \cite{FSV} and \cite{TV7}\).
\end{remark}
\begin{lemma}
\label{I-as}
Let \>$\Re z\to-\)\infty$ and $\mu$ is fixed. Then
\vvn.5>
\begin{align*}
& I_{a,b}(z,\mu\);m_1,m_2,l_1,l_2)\,=\,l_2!\>(l_2-b)!\;b!\,
(-\)\pi)^{-\)l_2}\exp\)(\mu\)z\)b/\ka)\,\x{}
\\[5pt]
&\>{}\x\,(-z/\ka)^{-\)(2\)b^2\]+\>b\)(m_1-\)m_2\)-\)2\)l_2)\)+\)m_2l_2)/\ka}
A_{l_2-b}(\mu\); m_1)\,A_b(\mu\);m_2)\bigl(\dl_{ab}+O(z^{-1})\bigr)\,.
\end{align*}
\end{lemma}
\nt
The lemma follows from \cite{TV1}\).

\subsection{Hypergeometric solutions of the trigonometric \KZ/ and difference
dynamical equations}
For $b=0\lc\min\)(m_2\),l_2)$, define
\vvn.4>
an \$M_{l_1}\]\ox L_{l_2}[m_1\),m_2]$-valued function
$$
\Jb_b(z,\mu\);l_1,l_2,m_1,m_2)\,=\>
\sum_{a=0}^{m_2}\,J_{a,b}(z,\mu\); l_1,l_2,m_1,m_2)
\;\frac1{(m_2-a)!\,a\)!}\,E^{m_2-a}_{2,1}\)v_{l_1} \ox E^{a}_{2,1}\)v_{l_2}\,.
\vv.3>
$$
Here the actual range of summation is until $a=\min\)(m_2,l_2)$,
since $E^{a}_{2,1}\)v_{l_2}=0$ for $a>l_2$.

\begin{theorem}
\label{J-KZ}
For any $b=0\lc\min\)(m_2,l_2)$ the function
\vvn.6>
\begin{align*}
U_b(z_1,z_2,\la_1,\la_2)\,=\,z_1^{(\)\la_1(l_1\]-\>m_2)\)+\)\la_2\)m_2\)-
\>m_2^2\)-\>m_1\)l_1\]+\>l_1^2/2)/\ka}\>
z_2^{\>l_2(\la_1-\>m_1+\>l_2/2)/\ka}\,\x\!\}{} &
\\[6pt]
{}\x\,(z_1-z_2)^{l_1l_2/\ka}\>\Jb_b(z_2/z_1,\la_2\)-\la_1;l_1,l_2,m_1,m_2) &
\\[-12pt]
\end{align*}
is a solution of the \KZ/ equations \Ref{KZ} and difference dynamical
equations \Ref{q-dyn} with values in $(M_{l_1}\]\ox L_{l_2})[m_1\),m_2]$.
\end{theorem}
\nt
The theorem is a direct corollary of \cite{MV}\).
\vsk.5>
We describe asymptotics of the integral $J_{a,b}(\emu,z\);l_1,l_2,m_1,m_2)$
as $\Re z\to-\)\infty$ and $\mu$ is fixed. For a positive number $\ka$,
a complex number $l$ and a nonnegative integer $m$ consider
the Selberg-type integral
\vvn.4>
$$
B_m(l)\,=\,
\int_{\gm_m}\,e^{-\!\sum_{u=1}^m\!s_u/\ka}\,\prod_{u=1}^m(-\)s_u)^{-1-\)l/\ka}
\prod_{1\le u<v\le m}(s_u\]-s_v)^{2/\ka}\,ds^{\)m}\,.
\vv.3>
$$
The integration contour $\gm_m$ has the form
\vvn.3>
$$
\gm_m\>=\>\lb\,(s_1\lc s_m)\in\C^{\>m}\vert s_u\in\Cc_u\,,
\ u=1\lc m\,\rb\,,
\vv.2>
$$
see the picture.

\vbox{\begin{center}
\begin{picture}(224,35)
\put(40,15){\oval(10,10)[l]}
\put(40,15){\circle*{3}}
\put(40,10){\vector(1,0){152}}
\put(40,20){\line(1,0){152}}
\qbezier[16](46,5)(118,5)(190,5)
\qbezier[16](46,25)(118,25)(190,25)
\qbezier[5](30,15)(30,25)(40,25)
\qbezier[4](40,5)(30,5)(30,15)
\put(40,15){\oval(30,30)[l]}
\put(40,0){\vector(1,0){152}}
\put(40,30){\line(1,0){152}}
\put(15,5){$0$}
\put(197,30){$s_m$}
\put(197,20){$s_1$}
\end{picture}
\vsk>
Picture 2. The contour $\gm_m$.
\end{center}}
\vsk.3>

\noindent
Here $\Cc_u,\, u=1\lc m$, are non-intersecting oriented loops in $\C$.
The loops start at $+\)\infty$, go around 0, and return to $+\)\infty$.
For $u<v$ the loop $\Cc_u$ lies inside the loop $\Cc_v$. We fix a univalued
branch of the integrand by assuming that at the point of $\gm_m$ where all
numbers $s_1\lc s_m$ are negative we have $\arg\)(-\)s_u)=0$ for
$u=1\lc m$, and $\arg\)(s_u\]-s_v)=0$ for $1\le u<v\le m$.
\vsk.2>
The formula for $B_m(l)$ is well known:
\begin{equation}
\label{selberg}
B_m(l)\,=\,(-2\pi i)^m\>\ka^{m\)(m\)-1-\>l)/\ka}\;
\prod_{j=0}^{m-1}\;\frac{\Gm(1-1/\ka)}{\Gm(1+(l-j)/\ka)\,\Gm(1-(j+1)/\ka)}\;,
\end{equation}
for example, \cf. \cite{TV2}\), \cite{MTV}\).

\begin{lemma}
\label{J-as}
Let \>$\Re z\to-\)\infty$ and $\mu$ is fixed. Then
\vvn.5>
\begin{align*}
J_{a,b}(\emu,z\); {}& l_1,l_2,m_1,m_2)\,=\,
(m_2-b)!\;b!\,e^{\)\pi\)i\)(m_2\)-\)b)\)(2b\>-\)l_2)/\ka}\,\x{}
\\[6pt]
& {}\!\!\}\x\,
e^{\>\mu\)b\)(z\)+\>l_1-\)2\)m_2\)+\)b)/\ka}\,
(1-\emu)^{-\)(2\)b^2\]+\>b\)(l_1-\>l_2\)-\)2\)m_2)\)+\)m_2\)l_2)/\ka}
\x{}
\\[5pt]
&{}\!\!\}\x\>(-z)^{-\)(2\)b^2\]+\>b\)(l_1-\>l_2\)-\)2\)m_2)\)+
\)m_2(m_2-\>l_1-1))/\ka}
B_{m_2-b}(\)l_1)\,B_b(\)l_2)\bigl(\dl_{ab}+O(z^{-1})\bigr)\,.\kern-.5em
\end{align*}
\end{lemma}
\vsk.2>
\nt
The proof is straightforward.

\subsection{Proof of Theorem \ref{first}}
Theorems~\ref{I-qKZ} and \ref{J-KZ}, and Corollary \ref{Xz} imply that for any
$b=0\lc\min\)(m_2,l_2)$ the functions $\Ib_b(z,\mu\);m_1,m_2,l_1,l_2)$ and
$X(z\);m_1,m_2)\>\phi^{-1}\bigl(\Jb_b(\emu,z\);l_1,l_2,\alb m_1, m_2)\bigr)$
satisfy the same first order difference equation with respect to $z$ with
step $\ka$. Hence, for any $a=0\lc\min\)(m_2,l_2)$ one has
\vv.3>
\begin{align}
\label{JIG}
X(z\);m_1,m_2)\> & J_{a,b}(\emu,z\);l_1,l_2,m_1,m_2)\,={}
\\[6pt]
{}=\!\sum_{c=0}^{\min\)(m_2,\>l_2)}\!{}
& I_{a,c}(z,\mu\);m_1,m_2,l_1,l_2)\>G_{b,c}(z,\mu\);m_1,m_2,l_1,l_2)\,,\!\!\!
\notag
\\[-13pt]
\notag
\end{align}
the connection coefficients $G_{b,c}(z,\mu\);m_1,m_2,l_1,l_2)$ being
\$\ka$-periodic functions of $z$ and holomorphic functions of $\mu$
in the strip $0<\Im\mu<2\pi$. Taking into account asymptotics of the integrals
$I_{a,c}(z,\mu\);m_1,m_2,l_1,l_2)$ and $J_{a,b}(\emu,z\);l_1,l_2,m_1,m_2)$ as
$\Re z\to-\)\infty$ and $\mu$ is fixed, see Lemmas~\ref{I-as} and \ref{J-as},
one can compute the connection coefficients and obtain formula \Ref{I=J}\).
Theorem~\ref{first} is proved.
\qed
\Rem
One can see from formula \Ref{I=J} that all connection coefficients
$G_{b,c}(z,\mu\);m_1,m_2,\alb l_1,l_2)$ in \Ref{JIG} as functions of $\mu$
are proportional to the same function $Y(\mu\);m_1,m_2,l_1,l_2)$. This fact,
which is pure computational in the given proof of Theorem~\ref{first},
can be observed independently in advance, because Theorems~\ref{J-KZ},
\ref{I-DD}, and Corollary \ref{KZDE} imply that the functions
$\Ib_b(z,\mu\);m_1,m_2,l_1,l_2)$ and $Y(\mu\);m_1,m_2,l_1,l_2)\>
\phi^{-1}\bigl(\Jb_b(\emu,z\);l_1,l_2,m_1,m_2)\bigr)$ satisfy
the same first order differential equation with respect to $\mu$.
\enddemo

\section{Differential dynamical operators}

Introduce the \em{dynamical differential operators} $D_1\>,\)D_2$ with
coefficients in $U(\glt)^{\ox\)n}\}$ acting on functions of complex variables
$z_1\lc z_n,\>\la_1,\la_2$ by the formula
\vvn.5>
\begin{align*}
& \,D_i(z_1\lc z_n,\la_1,\la_2)\,={}
\\[6pt]
&{}=\,\ka\>\la_i\frac\partial{\partial \la_i\]}\>+
\>\frac{\Et_{ii}^{\)2}}2\>-\)\tsum_{a=1}^n z_a\>E_{ii}^{\)(a)}\>-\)
\sum_{j=1}^2\,\sum_{1\le a<b\le n}\!E_{ij}^{\)(a)}\)E_{ji}^{\)(b)}\>-\)
\frac{\la_{i'}}{\la_i-\la_{i'}\]}\,(\Et_{21}\)\Et_{12}\)-\Et_{22})\,.
\\[-15pt]
\end{align*}
Here \>$\Et_{kl}=\sum_{a=1}^n\>E_{kl}^{\)(a)}$, and $i'$ is supplementary
to $i$, that is, $\lb i\>,i'\rb=\lb 1\>,2\rb$.
\begin{proposition}[\)\cite{TV5}\)]
One has $[\)Z_a\>,D_i\)]\)=\)0$ for all $a=1\lc n$
\vvn.1>
and $i=1,2$, where $Z_1\lc Z_n$ are the \qKZ/ operators \Ref{qKZo}\).
\end{proposition}
The differential equations
\vvn.4>
\begin{equation}
\label{dyn}
D_i(z_1\lc z_n,\la_1,\la_2)\,U(z_1\lc z_n,\la_1,\la_2)\,=\,0\,,\qquad i=1,2\,,
\kern-1em
\end{equation}
\vsk.2>\nt
for a \$V_1\lox V_n\)$-valued function $U(z_1\lc z_n,\la_1,\la_2)$ are called
the \em{dynamical differential equations} associated with the \qKZ/ equations.
\vsk.2>

\begin{theorem}[\)\cite{TV3}\)]
The isomorphism $\phi$, see~\Ref{phiso}, transforms the dynamical differential
operators acting in $(M_{m_1}\ox L_{m_2})[\)l_1\),l_2]$ into the trigonometric
\KZ/ operators \Ref{KZo} acting in $(M_{l_1}\ox L_{l_2})[m_1,m_2]${\rm:}
\vvn.2>
$$
\phi\;D_a(z_1,z_2,\la_1,\la_2)\,=\,
\nabla_a(\la_1,\la_2,z_1,z_2)\;\phi\,,\qquad a=1,2\,.
\vvgood
\vv.4>
$$
\end{theorem}

\begin{corollary}
\label{KZDE}
Let an \$(M_{l_1}\]\ox L_{l_2})[m_1,m_2]$-valued function
$U(z_1,z_2,\la_1,\la_2)$ solve the tri\-go\-no\-metric \KZ/ equations\/{\rm:}
\vvn.1>
$$
\nabla_a(z_1,z_2,\la_1,\la_2)\,U(z_1,z_2,\la_1,\la_2)\,=\,0\,, \qquad a=1,2\,.
\vv.4>
$$
Then the \$(M_{m_1}\]\ox L_{m_2})[\)l_1\),l_2]$-valued function
\vvn.4>
$$
\Ut(z_1,z_2,\la_1,\la_2)\,=\,\phi^{-1}\bigl(U(\la_1,\la_2,z_1,z_2)\bigr)
\vv.4>
$$
solves the system of the dynamical differential equations \Ref{dyn}\){\rm:}
\vvn.4>
$$
D_a(z_1,z_2,\la_1,\la_2)\,\Ut(z_1,z_2,\la_1,\la_2)\,=\,0\,, \qquad a=1,2\,.
\vv.2>
$$
\end{corollary}

The next statement describes \q-hypergeometric solutions of the dynamical
differential equations.
\begin{theorem}
\label{I-DD}
For any $b=0\lc\min\)(m_2,l_2)$ the function $U_b(z_1,z_2,\la_1,\la_2)$
defined by \Ref{Ub} is a solution of equations \Ref{dyn} with values in
$(M_{m_1}\]\ox L_{m_2})[\)l_1\),l_2]$.
\end{theorem}
\nt
The theorem follows from \cite{TV6}\).

\end{document}